\documentclass[12TP,draft]{article}

\begin{document}

\begin{center}
\LARGE\noindent\textbf{ On Hamiltonian Bypasses in Digraphs with the Condition of Y. Manoussakis}\\

\end{center}
\begin{center}
\noindent\textbf{Samvel Kh. Darbinyan}\\

Institute for Informatics and Automation Problems, Armenian National Academy of Sciences

E-mail: samdarbin@ipia.sci.am \\
\end{center}

\textbf{Abstract}

Let $D$ be a strongly connected directed graph of order $n\geq 4$ vertices which satisfies the following condition for every triple $x,y,z$ of vertices such that $x$ and $y$ are non-adjacent: If there is no arc from $x$ to $z$, then $d(x)+d(y)+d^+(x)+d^-(z)\geq 3n-2$. If there is no arc from $z$ to $x$, then $d(x)+d(y)+d^-(x)+d^+(z)\geq 3n-2$. In \cite{[15]} (J. of Graph Theory, Vol.16, No. 5,  51-59, 1992) Y. Manoussakis proved that $D$ is Hamiltonian. In [9] it was shown that $D$ contains a pre-Hamiltonian cycle (i.e., a cycle of length $n-1$) or $n$ is even and $D$ is isomorphic to the complete bipartite digraph with partite sets of cardinalities of $n/2$ and $n/2$. In this paper we show that $D$ contains also a Hamiltonian bypass, (i.e., a subdigraph  obtained from a Hamiltonian cycle by reversing exactly one arc) or $D$ is isomorphic to one tournament of order 5.\\ 

\textbf{Keywords:} Digraphs, cycles, Hamiltonian cycles, Hamiltonian bypasses. \\

\section {Introduction} 

The directed graph (digraph) $D$ is Hamiltonian  if it contains a Hamiltonian cycle,
i.e., a cycle that includes every vertex of $D$. A Hamiltonian bypass in $D$ is a subdigraph obtained from a Hamiltonian cycle by reversing exactly one arc. We recall the following well-known degree conditions (Theorems 1-6) that guarantee that a digraph is Hamiltonian. 

 \textbf{Theorem 1} (Nash-Williams \cite{[17]}). Let $D$ be a digraph of order $n$ such that for every vertex $x$, $d^+(x)\geq n/2$ and $d^-(x)\geq n/2$, then $D$ is Hamiltonian.\\

 \textbf{Theorem 2} (Ghouila-Houri \cite{[14]}). Let $D$ be a strong digraph of order $n$. If $d(x)\geq n$ for all vertices $x\in V(D)$, then $D$ is Hamiltonian.\\

 \textbf{Theorem 3} (Woodall \cite{[19]}). Let $D$ be a digraph of order $n\geq 2$. If $d^+(x)+d^-(y)\geq n$ for all pairs of vertices $x$ and $y$ such that there is no arc from $x$ to $y$, then $D$ is Hamiltonian.\\
 
\textbf{Theorem 4} (Meyniel \cite{[16]}). Let $D$ be a strong digraph of order $n\geq 2$. If $d(x)+d(y)\geq 2n-1$ for all pairs of non-adjacent vertices in $D$, then $D$ is Hamiltonian.

It is easy to see that Meyniel's theorem is a common generalization of Ghouila-Houri's and Woodall's theorems. For a short proof of Theorem 1.3, see \cite{[5]}. \\

C. Thomassen \cite{[18]} (for $n=2k+1$) and S. Darbinyan \cite{[6]}
(for $n=2k$) proved the following:

\textbf{Theorem 5} \cite{[18], [6]}. If $D$ is a digraph of order
$n\geq 5$ with minimum degree at least $n-1$ and with minimum
semi-degree at least $n/2-1$, then $D$ is Hamiltonian (unless some
extremal cases which are characterized).\\

In view of the next theorems we need the following definitions.

\textbf{Definition 1} \cite{[15]}. Let $k$ be an integer. A digraph $D$ of order $n\geq 3$ satisfies the condition $A_k$ if and only if for every triple of vertices $x,y,z$ such that $x$ and $y$ are non-adjacent: If there is no arc from $x$ to $z$, then $d(x)+d(y)+d^+(x)+d^-(z)\geq 3n-2+k$. If there is no arc from $z$ to $x$, then $d(x)+d(y)+d^-(x)+d^+(z)\geq 3n-2+k$.  \\

\textbf{Definition 2}. Let $D_0$ denote any digraph of order $n\geq 5$, $n$ odd, such that $V(D_0)=A\cup B$, where $A\cap B=\emptyset$, $A$ is an independent set with $(n+1)/2$ vertices, $B$ is a set of $(n-1)/2$ vertices inducing any arbitrary subdigraph, and $e(A,B)=(n+1)(n-1)/2$. 
$D_0$ satisfies the condition $A_{-1}$, but has no Hamiltonian bypass. \\

\textbf{Definition 3}. For any $k\in [1, n-2]$ let $D_1$ denote a digraph of order $n\geq 4$, obtained from $K^*_{n-k}$ and $K^*_{k+1}$ by identifying a vertex of the first with a vertex of the second. 
$D_1$ satisfies the condition $A_{-1}$, but has no Hamiltonian bypass.\\

\textbf{Definition 4}. By $T(5)$ we denote a tournament of order 5 with vertex set $V(T(5))=\{x_1, x_2, x_3, x_4, y\}$ and arc set $A(T(5))=\{x_ix_{i+1}/ i\in [1,3]\}\cup \{x_4x_1, x_1y, x_3y, yx_2, yx_4,x_1x_3, x_2x_4\}$. $T(5)$ satisfies condition $A_0$, but has no Hamiltonian bypass.\\

\textbf{Theorem 6} (Manoussakis \cite{[15]}). If a strong digraph $D$  satisfies the condition $A_0$, then $D$ is Hamiltonian. \\

In \cite{[4]} it was proved that if a digraph $D$ satisfies the condition of Nash-Williams' or Ghouila-Houri's or Woodall's theorem, then $D$ contains a Hamiltonian bypass. In [4] also proved the following theorem:

 \textbf{Theorem 7} (Benhocine \cite{[4]}). Every strongly 2-connected digraph of order $n$ and minimum degree at least $n-1$ contains a Hamiltonian bypass, unless $D$ is isomorphic to a digraph of type $D_0$.\\

In \cite{[7]} the following theorem was proved:

 \textbf{Theorem 8} (Darbinyan \cite{[7]}). Let $D$ be a strong digraph of order $n\geq 3$. If $d(x)+d(y)\geq 2n-2$ for all pairs of non-adjacent vertices in $D$, then $D$ contains a Hamiltonian bypass unless it is isomorphic to a digraph of the set $D_0\cup \{D_1, T_5, C_3\}$, where $C_3$ is a directed cycle of length 3.\\

For $n\geq 3$ and $k\in [2,n]$, $D(n,k)$ denotes the digraph of order $n$ obtained from a directed cycle $C$ of length $n$ by reversing exactly $k-1$ consecutive arcs. In \cite{[7],[8]} Darbinyan  studied the problem of the existence of $D(n,3)$ in digraphs with condition of Meyniel's theorem and in oriented graphs with large in-degrees and out-degrees.

 \textbf{Theorem 9} (Darbinyan \cite{[7]}). Let $D$ be a strong digraph of order $n\geq 4$. If $d(x)+d(y)\geq 2n-1$ for all pairs of non-adjacent vertices in $D$, then $D$ contains a $D(n,3)$.\\

\textbf{Theorem 10} (Darbinyan \cite{[8]}). Let $D$ be a oriented graph of order $n\geq 10$. If the minimum in-degree and out-degree of $D$ at least $(n-3)/2$, then $D$ contains a $D(n,3)$.\\

In \cite{[9]} the following theorem was proved: 

 \textbf{Theorem 11}. Any strongly connected digraph $D$ of order $n\geq 4$ satisfying the condition $A_0$ contains a pre-Hamiltonian cycle (i.e., a cycle of length $n-1$) or $n$ is even and $D$ is isomorphic to the complete bipartite digraph with partite sets of cardinalities $n/2$ and $n/2$. \\

In this paper using Theorem 11 we prove the following:\\

 \textbf{Theorem 12}. Any strongly connected digraph $D$ of order $n\geq 4$ satisfying the condition $A_0$ contains a Hamiltonian bypass unless  $D$ is isomorphic to the tournament $T(5)$. \\

The following two examples show the sharpness of the condition of Theorem 10. The digraph consisting of the disjoint union of two complete digraphs with one common vertex shows that the bound in the above theorem is best possible and the digraph obtained from a complete bipartite digraph after deleting one arc.\\

\section {Terminology and Notations}

We shall assume that the reader is familiar with the standard
terminology on the directed graphs (digraph)
 and refer the reader to the  monograph of Bang-Jensen and Gutin \cite{[1]} for terminology not discussed here.
  In this paper we consider finite digraphs without loops and multiple arcs. For a digraph $D$, we denote
  by $V(D)$ the vertex set of $D$ and by  $A(D)$ the set of arcs in $D$. The order of $D$ is the number
  of its vertices. Often we will write $D$ instead of $A(D)$ and $V(D)$. The arc of a digraph $D$ directed from
   $x$ to $y$ is denoted by $xy$. For disjoint subsets $A$ and  $B$ of $V(D)$  we define $A(A\rightarrow B)$ \,
   as the set $\{xy\in A(D) / x\in A, y\in B\}$ and $A(A,B)=A(A\rightarrow B)\cup A(B\rightarrow A)$. If $x\in V(D)$
   and $A=\{x\}$ we write $x$ instead of $\{x\}$.  The out-neighborhood of a vertex $x$ is the set $N^+(x)=\{y\in V(D) / xy\in A(D)\}$ and $N^-(x)=\{y\in V(D) / yx\in A(D)\}$ is the in-neighborhood of $x$. Similarly, if $A\subseteq V(D)$, then $N^+(x,A)=\{y\in A / xy\in A(D)\}$ and $N^-(x,A)=\{y\in A / yx\in A(D)\}$. The out-degree of $x$ is $d^+(x)=|N^+(x)|$ and $d^-(x)=|N^-(x)|$ is the in-degree of $x$. Similarly, $d^+(x,A)=|N^+(x,A)|$ and $d^-(x,A)=|N^-(x,A)|$. The degree of the vertex $x$ in $D$ defined as $d(x)=d^+(x)+d^-(x)$ (similarly, $d(x,A)=d^+(x,A)+d^-(x,A)$). The subdigraph of $D$ induced by a subset $A$ of $V(D)$ is denoted by $\langle A\rangle$. The path (respectively, the cycle) consisting of the distinct vertices $x_1,x_2,\ldots ,x_m$ ( $m\geq 2 $) and the arcs $x_ix_{i+1}$, $i\in [1,m-1]$  (respectively, $x_ix_{i+1}$, $i\in [1,m-1]$, and $x_mx_1$), is denoted  $x_1x_2\cdots x_m$ (respectively, $x_1x_2\cdots x_mx_1$). We say that $x_1x_2\cdots x_m$ is a path from $x_1$ to $x_m$ or is an $(x_1,x_m)$-path. For a cycle  $C_k:=x_1x_2\cdots x_kx_1$ of length $k$, the subscripts considered modulo $k$, i.e. $x_i=x_s$ for every $s$ and $i$ such that  $i\equiv s\, (\hbox {mod} \,k)$. A cycle that contains the all vertices of $D$ (respectively, the all vertices of $D$ except one) is a Hamiltonian cycle (respectively, is a pre-Hamiltonian cycle). If $P$ is a path containing a subpath from $x$ to $y$ we let $P[x,y]$ denote that subpath. Similarly, if $C$ is a cycle containing vertices $x$ and $y$, $C[x,y]$ denotes the subpath of $C$ from $x$ to $y$. A digraph $D$ is strongly connected (or, just, strong) if there exists a path from $x$ to $y$ and a path from $y$ to $x$ for every pair of distinct vertices $x,y$.
  For an undirected graph $G$, we denote by $G^*$ symmetric digraph obtained from $G$ by replacing every edge $xy$ with the pair $xy$, $yx$ of arcs.  $K_{p,q}$ denotes the complete bipartite graph  with partite sets of cardinalities $p$ and $q$.  Two distinct vertices $x$ and $y$ are adjacent if $xy\in A(D)$ or $yx\in A(D) $ (or both). For integers $a$ and $b$, $a\leq b$, let $[a,b]$  denote the set of
all integers which are not less than $a$ and are not greater than
$b$.  By $D(n;2)=[x_1x_{n}; x_1x_2\ldots , x_n]$ is denoted the Hamiltonian bypass obtained from a Hamiltonian cycle $x_1x_2\ldots x_nx_1$ by reversing the arc $x_nx_1$.

\section { Preliminaries }

The following well-known simple Lemmas 1-4 are the basis of our results
and other theorems on directed cycles and paths in digraphs. They
will be used extensively in the proof of our result. \\

\textbf{Lemma 1} \cite{[13]}. Let $D$ be a digraph of order  $n\geq 3$
 containing a
 cycle $C_m$, $m\in [2,n-1] $. Let $x$ be a vertex not contained in this cycle. If $d(x,C_m)\geq m+1$,
 then  $D$ contains a cycle $C_k$ for all  $k\in [2,m+1]$.\\

The following lemma is a slight modification of a lemma by Bondy and Thomassen \cite{[5]}.

\textbf{Lemma 2}. Let $D$ be a digraph of order $n\geq 3$
 containing a
 path $P:=x_1x_2\ldots x_m$, $m\in [2,n-1]$ and let $x$ be a vertex not contained in this path.
  If one of the following conditions holds:

 (i) $d(x,P)\geq m+2$;

 (ii) $d(x,P)\geq m+1$ and $xx_1\notin D$ or $x_mx_1\notin D$;

 (iii) $d(x,P)\geq m$, $xx_1\notin D$ and $x_mx\notin D$,

then there is an  $i\in [1,m-1]$ such that
$x_ix,xx_{i+1}\in D$ (the arc $x_ix_{i+1}$ is a partner of $x$), i.e., $D$ contains a path $x_1x_2\ldots
x_ixx_{i+1}\ldots x_m$ of length $m$ (we say that  $x$ can be
inserted into $P$ or the path $x_1x_2\ldots x_ixx_{i+1}\ldots$ $ x_m$
is extended from $P$ with $x$ ).\\

If in Lemma 1 and  Lemma 2 instead of the vertex $x$ consider a path $Q$, then we get the following Lemmas 3 and 4, respectively. \\

\textbf{Lemma 3}. Let $C_k:=x_1x_2\ldots x_kx_1$, $k\geq 2$, be a non-Hamiltonian cycle in a digraph $D$. Moreover, assume that there exists a path $Q:=y_1y_2\ldots y_r$, $r\geq 1$, in $D-C_k$. If $d^-(y_1,C_k)+d^+(y_r,C_k)\geq k+1$, then for all $m\in [r+1,k+r]$ the digraph $D$ contains a cycle $C_m$ of length $m$  with vertex set $V(C_m)\subseteq V(C_k)\cup V(Q)$.\\

\textbf{Lemma 4}. Let $P:=x_1x_2\ldots x_k$, $k\geq 2$, be a non-Hamiltonian path in a digraph $D$. Moreover, assume that there exists a path $Q:=y_1y_2\ldots y_r$, $r\geq 1$, in $D-P$. If $d^-(y_1,P)+d^+(y_r,P)\geq k+|A(y_1\rightarrow x_1)|+|A(x_k\rightarrow y_r)|$, then there is a $x_i$, $i\in [1,k-1]$, such that $x_iy_1, y_rx_{i+1}\in D$ and $D$ contains a path from $x_1$ to $x_k$ with vertex set $V(P)\cup V(Q)$.\\

In the proof of Theorem 11 we also need the following
 lemma which is a simple extension of a lemma by Y. Manoussakis \cite{[15]}.

\textbf{Lemma 5}. Let $D$ be a digraph of order $n\geq 3$ satisfying condition $A_0$. Assume that there are two distinct pairs $x,y$ and $x,z$ of non-adjacent vertices in $D$. If $d(x)+d(y)\leq 2n-a$ for some integer $a\geq 1$, then  $d(x)+d(z)\geq 2n-2+a/2$. In particular, if $d(x)+d(y)\leq 2n-2$, then $d(x)+d(z)\geq 2n-1$.\\ 

\textbf{Definition 5} (\cite{[1]}, \cite{[2]}). Let $Q=y_1y_2\ldots y_s$ be a path in a digraph $D$ (possibly, $s=1$) and let $P=x_1x_2\ldots x_t$, $t\geq 2$, be a path in $D-V(Q)$. $Q$ has a partner on $P$ if there is an arc (the partner of $Q$) $x_ix_{i+1}$ such that $x_iy_1,y_sx_{i+1}\in D$. In this case the path $Q$ can be inserted into $P$ to give a new $(x_1,x_t)$-path with  vertex set $V(P)\cup V(Q)$. The path $Q$ has a collection of partners on $P$ if there are integers $i_1=1<i_2<\cdots <i_m=s+1$ such that, for every $k=2,3,\ldots , m$ the subpath $Q[y_{i_{k-1}},y_{i_k-1}]$ has a partner on $P$.\\

\textbf{Lemma 6} (\cite{[1]}, \cite{[2]}, Multi-Insertion Lemma).  Let $Q=y_1y_2\ldots y_s$ be a path in a digraph $D$ (possibly, $s=1$) and let $P=x_1x_2\ldots x_t$, $t\geq 2$, be a path in $D-V(Q)$. If $Q$ has a collection of partners on $P$, then there is an $(x_1,x_t)$-path with vertex set $V(P)\cup V(Q)$.\\

The following lemma is obvious.

\textbf{Lemma 7}. Let $D$ be a digraph of order $n\geq 3$ and let $C:=x_1x_2\ldots x_{n-1}x_1$ be an arbitrary cycle of length $n-1$ in $D$. If a vertex $y$ is not on $C$ and $D$ contains no Hamiltonian bypass, then

(i) $d^+(y,\{x_i, x_{i+1}\})\leq 1$ and   $d^-(y,\{x_i, x_{i+1}\})\leq 1$ for all $i\in [1,n-1]$;

(ii) $d^+(y)\leq (n-1)/2$, $d^-(y)\leq (n-1)/2$ and $d(y)\leq n-1$;

(iii) if $x_ky, yx_{k+1}\in D$, then $x_{i+1}x_i\notin D$ for all $x_i\not=x_k$.\\

Let $D$ be a digraph of order $n\geq 3$ and let $C_{n-1}$ be a cycle of length $n-1$  
in $D$. If for the vertex $y\notin C_{n-1}$, $d(y)\geq n$, then we say that $C_{n-1}$ is a good cycle. Notice that, by Lemma 7, if a digraph $D$ contains a good cycle, then $D$ also contains a Hamiltonian bypass.
 
\section {Proof of Theorem 12}

In the proof of Theorem 12 we often will use the following definition:

\textbf{Definition 6}. Let $P_0:=x_1x_2\ldots x_m$, $m\geq 2$, be an $(x_1,x_m)$-path in $D$ and let the vertices
$y_1,y_2,\ldots y_k\in V(D)- V(P_0)$. For $i\in [1,k]$ we denote by $P_i$ an $(x_1,x_m)$-path in $D$ with vertex set $V(P_{i-1})\cup \{y_j\}$ (if it exists), i.e, $P_i$ is extended path obtained from $P_{i-1}$ with some vertex $y_j$, where $y_j\notin V(P_{i-1}$). If $e+1$ is the maximum possible number of these paths $P_0, P_1,\ldots , P_e$, $e\in [0,k]$, then we say that $P_e$ is extended path obtained from $P_0$ with vertices $y_1,y_2,\ldots , y_k$ as much as possible. Notice that $P_i$ is an $(x_1,x_m)$-path of length $m+i-1$ for all $i\in [0,e]$.\\

\textbf{Proof of Theorem 12}. By Theorem 9 the digraph $D$ contains a cycle of length $n-1$ or $n$ is even and $D$ is isomorphic to the complete bipartite digraph with partite sets of cardinalities of $n/2$ and $n/2$. If $D$ is a complete bipartite digraph then it is easy to see that $D$ has a Hamiltonian bypass.  In the sequel, we assume that $D$ contains a cycle of length $n-1$. Let $C=x_1x_2\ldots x_{n-1}x_1$ be an arbitrary cycle of length $n-1$ in $D$ and let $y\notin C$. It is a simple matter to check that for $n=4$ the theorem is true. Further, let $n\geq 5$. Note that from the condition $A_0$ and Lemma 5 immediately follows that $d(y)\geq 3$. Now suppose, to the contrary, that $D$ contains no Hamiltonian bypass (by Lemma 7(ii) it is clear that $D$ also contains no good cycle).

For the cycle $C$ and the vertex $y$ we prove the following Claims 1- 7 below.\\

 \textbf{Claim 1}. $d(y,\{x_i\})\leq 1$ for all $i\in [1,n-1]$.

\textbf{Proof}. Assume that the claim is not true. Without loss of generality, assume that $d(y,\{x_{n-1}\})=2$, i.e., $x_{n-1}y, yx_{n-1}\in D$. By Lemma 7(i), $y$  is not adjacent with $x_1$ and $x_{n-2}$. Since $d(y)\geq 3$, we can assume that for some integers $a\geq 1$ and $b\geq 1$ the following hold
$$ d(y, \{x_1,x_2, \ldots , x_a\})=d(y,\{x_{n-2},x_{n-3}, \ldots ,x_{n-b-1}\})=0, \eqno (1) $$
and
$$ min\{d(y,\{x_{a+1}\}), d(y,\{x_{n-b-2}\})\} \geq 1 \eqno (2) 
$$
($x_{n-b-2}=x_{a+1}$ is possible). Now from  Lemma 7(i) and (1) it follows that 
$$d(y)=d(y,\{x_{n-1}\})+ d(y,C[x_{a+1},x_{n-b-2}])\leq n-b-a+1.  \eqno (3) $$
If there is an $(x_{a+1},x_{n-1})$-path $P$ (respectively, an $(x_{n-1}, x_{n-b-2}$)-path $Q$) with vertex set $V(C)$, then, since (2) and $d(y,\{x_{n-1}\})=2$, it is easy to see that $D$ contains a Hamiltonian bypass. So we may assume that there is no $(x_{a+1},x_{n-1})$-path and there is no $(x_{n-1}, x_{n-b-2}$)-path with vertex set $V(C)$. We extend the path $P_0:=C[x_{a+1},x_{n-1}]$ (respectively, $P_0:=C[x_{n-1},x_{n-b-2}]$) with vertices $ x_1,x_2,\ldots , x_a$ (respectively, $x_{n-b-1},x_{n-b},\ldots , x_{n-2}$) as much as possible. Then some vertices $z_1, z_2, \ldots , z_d\in \{x_1,x_2,\ldots , x_a\}$, $d\in [1,a]$, (respectively, $u_1, u_2, \ldots , u_l\in \{x_{n-b-1},x_{n-b},\ldots , x_{n-2}\}$, $l\in [1,b]$) are not on the extended path $P_e$. Therefore using Lemma 2(i), we obtain that 
$$ 
d(z_i)\leq n+d-2 \quad \hbox {and} \quad d(u_j)\leq n+l-2 \eqno (4) $$ 
for all $i\in [1,d]$ and $j\in [1,l]$.  Since $d\leq a+b-1$ and $l\leq a+b-1$, from inequalities (3) and (4) it follows that 
$$
d(y)+d(z_i)\leq 2n-1+d-a-b\leq 2n-2 \quad and \quad d(y)+d(u_j)\leq 2n-1+l-a-b\leq 2n-2.
$$
The last two inequalities contradicts Lemma 5. Claim 1 is proved.  \fbox \\\\

\textbf{Claim 2}. $d(y)\leq n-2$.

\textbf{Proof}. Suppose, on the contrary, that $d(y)\geq n-1$. Then, by Lemma 7(ii), $d(y)=n-1$. Using Lemma 7(i) and Claim 1, we obtain that $n$ odd ($n:=2m+1$), and without loss of generality, we may assume that
$$
N^+(y)=\{x_1,x_3,\ldots , x_{n-2}\} \quad \hbox{and} \quad N^-(y)=\{x_2,x_4,\ldots , x_{n-1}\}.   \eqno (5) $$

By Lemma 7(iii), 
$$x_{i+1}x_i\notin D \quad \hbox{for all} \quad i\in [1,n-1]. \eqno (6) $$

\textbf{Case 2.1}. There is a $x_i$ such that $d(x_i)\geq n$. Without loss of generality, we may assume that $d(x_1)\geq n$ because of (5). Since $D$ contains no Hamiltonian bypass, it follows that $x_1$ has no partner on $C[x_3, x_{n-2}]$. From (6), Lemma 2(ii) and 
$$ n\leq d(x_1)= d(x_1, \{x_2, x_{n-1}, y\})+ d(x_1,C[x_3, x_{n-2}])
$$
it follows that $d(x_1,C[x_3, x_{n-2}])= n-3$ and $x_1x_3, x_{n-2}x_1\in D$. If $x_{n-1}x_2\in D$, then $D(n,2)=[x_{n-1}x_2;$ $ x_{n-1}yx_3x_4\ldots x_{n-2}x_1x_2]$, and if $x_{2}x_{n-1}\in D$, then $D(n,2)=[x_{2}x_{n-1}; x_{2}yx_1x_3x_4\ldots x_{n-1}]$, which contradicts to our assumption. So, we can assume that $x_2, x_{n-1}$ are non-adjacent. Since $yx_1x_3x_4 \ldots x_{n-1}y$ (respectively, $x_{n-2}x_1x_2yx_3\ldots x_{n-2}$) is a cycle of length $n-1$ which does not contain the vertex $x_2$ (respectively, $x_{n-1}$), by Lemma 7(ii), $d(x_2)\leq n-1$  (respectively,  $d(x_{n-1})\leq n-1$) and $d^-(x_2)\leq (n-1)/2=m$. Now since the triple of vertices $x_{n-1}, x_2, y$ satisfies the condition $A_0$, we obtain that
$$
3n-2\leq d(x_{n-1})+d(x_2)+d^-(x_2)+d^+(y)\leq 2n-2+2m=3n-3,
$$
which is a contradiction.

\textbf{Case 2.2}. $d(x_{i})\leq n-1$ for all $i\in [1,n-1]$. Observe that $d(x_i)+d(x_j)\leq 2n-2$ for all distinct vertices $x_i$ and $x_j$. Observe that this together with Lemma 5 implies that every vertex $x_i$ is adjacent with all vertices of $D$ maybe except only one vertex.

 \textbf{Subcase 2.2.1}. $x_ix_{i+2}\in D$ for some $i\in [1,n-1]$. Without loss of generality, assume that $x_1x_3\in D$. Then 

(i) $x_2x_4\notin D$ (otherwise, if $x_2x_4\in D$, then $D(n,2)=[x_2x_3; x_2x_4x_5\ldots x_{n-1}yx_1x_3]$).

(ii) $x_2x_{n-1}\notin D$ (otherwise, if $x_2x_{n-1}\in D$, then $D(n,2)=[x_2x_{n-1}; x_2yx_1x_3\ldots x_{n-1}]$).

(iii) $x_{n-1}x_{2}\notin D$ (otherwise, if $x_{n-1}x_{2}\in D$ and $n\geq 6$, then $D(n,2)=[x_1x_{2}; x_1x_3x_4yx_5\ldots x_{n-1}x_2]$), and if $x_{n-1}x_2\in D$ and $n=5$, then $x_3x_1\notin D$ and $D$ is isomorphic to $T(5)$).

Therefore, if $D$ is not isomorphic to $T(5)$, then by (ii) and (iii), $x_2, x_{n-1}$ are non-adjacent. Now we will consider the cycle $C':=yx_1x_3x_4\ldots x_{n-1}y$ of length $n-1$ which doese not contain $x_2$. By Lemma 7(ii), $d^-(x_2)\leq m$. This together with $d^+(y)=m$, $d(x_2)$ and $d(x_{n-1})\leq n-1$ implies that 
$$
d(x_{n-1})+d(x_2)+d^-(x_2)+d^+(y)\leq 2n-2+2m=3n-3,
$$
which contradicts the condition $A_0$, since $x_2, x_{n-1}$ are non-adjacent and $yx_2\notin D$. 

  \textbf{Subcase 2.2.2}. $x_ix_{i+2}\notin D$ for all $i\in [1,n-1]$. It is not difficult to see that any $x_i$ cannot be inserted into $C[x_{i+1},x_{i-1}]$. By Lemma 2(iii), $d(x_i, C[x_{i+2},x_{i-2}])\leq n-5$. Therefore, since $d(x_i,\{y,x_{i-1}, x_{i+1}\})=3$, we have that $d(x_i)\leq n-2$ for all $i\in [1,n-1]$. By Lemma 5, from this and the above observation we conclude that $D$ contains no cycle of length two, every vertex $x_i$ is adjacent exactly with $n-2$ vertices, and hence $d(x_i)=n-2$ for all $x_i$.

First we consider the vertex $x_2$. Without loss of generality, assume that $x_2, x_r$ are non-adjacent, where $r\in [4, n-1]$. The triple of vertices $x_2,x_r, y$ satisfies the condition $A_0$, since $yx_2\notin D$. Therefore 

$$3n-2\leq d(x_{r})+d(x_2)+d^-(x_2)+d^+(y)\leq 2n-4+(n-1)/2+d^-(x_2) \eqno(7) $$
and $d^-(x_2)\geq (n+5)/2=m+3$ (recall that $n=2m+1$). From this, since $x_2$ cannot be inserted into $C[x_3,x_1]$ and $x_2x_4\notin D$, $x_{n-1}x_2\notin D$, we obtain that 
$$
N^-(x_2)=\{x_1,x_4,x_5,\ldots , x_{r-1}\} \quad \hbox{and} \quad N^+(x_2)=\{y,x_3,x_{r+1},x_{r+2},\ldots , x_{n-1}\}. \eqno (8)$$

In particular, $r\geq m+6$ and $x_4x_2\in D$. Now we consider the vertex $x_1$.  Without loss of generality, assume that $x_1, x_k$ are non-adjacent, where $k\in [3,n-2]$. Similarly (7) and (8), we obtain  
$$
3n-2\leq d(x_{1})+d(x_k)+d^+(x_1)+d^-(y),\quad d^+(x_1)\geq m+3, $$
$$
N^+(x_1)=\{x_2,x_{n-2},x_{n-3},\ldots , x_{k+1}\}\quad \hbox{and}\quad k\leq r-1 .$$
In particular, $x_1x_r\in D$. By symmetry of $x_1$ and $x_3$, we also have that $x_3x_{n-1}\in D$. Now from (5) and (8) we have that $D(n,2)=[x_3x_{n-1}; x_3x_4\ldots x_{r-1}x_2yx_1x_r\ldots x_{n-1}]$. This is contrary to the our assumption and completes the proof of Claim 2. \fbox \\\\

\textbf{Claim 3}. Let $d(y,C[x_{l+1},x_{k-1}])=0$ and $y$ is adjacent with $x_l$ and $x_k$, where $a+2:=|C[x_{l},x_{k}]|\geq 3$. Then

(i) if $x_ly, x_ky\in D$ or $yx_l, yx_k\in D$, then there is a vertex $u\in C[x_{l+1},x_{k-1}]$ such that $d(y)+d(u)\leq 2n-3$;

(ii) if $x_ly, yx_k\in D$, then there is an $(x_{k},x_{l})$-path with vertex set $V(C)-\{u\}$, where $u$ is some vertex of $C[x_{l+1},x_{k-1}]$ and   $d(u)\leq n-1$. In particular, $d(y)+d(u)\leq 2n-3$.

(iii) if $x_ly, yx_k\in D$ (or $yx_l, yx_k\in D$ or $x_ly, x_ky\in D$), then there are no $x_i$ and $x_j$ such that $C[x_i,x_j]\not= C[x_l,x_k]$, $b:=|C[x_i,x_j]|\geq 3$, $d(y, C[x_{i+1},x_{j-1}])=0 $ and a) $x_iy,x_jy\in D$ or b) $yx_i,yx_j\in D$ or c) $x_iy,yx_j\in D$.

\textbf{Proof}. By Claim 1, $d(y)\leq n-a-1$.

(i). It is not difficult to see that there is no $(x_k,x_l)$-path with vertex set $V(C)$. We extend the path $P_0:=C[x_k,x_l]$ with vertices $x_{l+1},x_{l+2},\ldots, ,x_{k-1}$ as much as possible. Then some vertices $z_1,z_2,\ldots , z_d\in \{x_{l+1},x_{l+2},\ldots , x_{k-1}\}$, $d\in [1,a]$, are not on the obtained extended path $P_e$. Hence using Lemma 2(i) we obtain that $d(z_i)\leq n+d-2$ (let $u:=z_1$). Therefore for all $i\in [1,d]$
$$d(y)+d(z_i)\leq n-a-1+n+d-2\leq 2n-3. \eqno (9) $$

(ii). Assume, without loss of generality, that $x_{n-1}y, yx_{a+1}\in D$ (i.e., $x_l=x_{n-1}$ and $x_k=x_{a+1}$) and  $d(y,C[x_{1},x_{a}])=0$ where $a\in [1,n-4]$. 
If $a=1$, then Claim 3(ii)  clearly is true. So, we can assume that $a\geq 2$. We extend the path $P_0:=C[x_{a+1},x_{n-1}]$ with vertices $x_1,x_2,\ldots ,x_a$ as much as possible. Then some vertices $z_1,z_2,\ldots , z_d\in \{x_1,x_2,\ldots ,x_a\}$ are not in the extended path $P_e$. We claim that $d=0$ or $d=1$. Indeed, if $d\geq 2$, then for the vertices $z_1$ and $z_2$ inequality (9) holds,  which contradicts Lemma 5. Therefore $d=0$ or $d=1$. 
If $d=1$, then $d(z_1)\leq n-1$ (let $u:=z_1$) and $P_e$ is an $(x_{a+1},x_{n-1})$-path with vertex set $V(C)-\{u\}$, and if $d=0$, then $e\geq 2$, $P_{e-1}$ is an  $(x_{a+1},x_{n-1})$-path with vertex set $V(C)- \{u \}$, where now $u$ is some vertex of $C[x_1,x_a]$, and  $d(u)\leq n-1$ since $D$ contains no good cycle. It is clear that  $d(y)+d(u)\leq 2n-3$.

(iii). Assume that Claim 3(iii) is not true. From Claims 3(i) and 3(ii) it follows that there are two distinct vertices  $u\in C[x_{l+1},x_{k-1}]$ and   $v\in C[x_{i+1},x_{j-1}]$ such that  $d(y)+d(u)\leq 2n-3$ and $d(y)+d(v)\leq 2n-3$. These last two inequalities contradicts Lemma 5, since $y,u$  and $y,v$ are two distinct pairs of non-adjacent vertices. Claim 3 is proved. \fbox \\\\

\textbf{Claim 4}. There are no two distinct vertices $x_i$ and $x_j$ such that $x_iy, x_jy\in D$ (or $yx_i, yx_j\in D$), $|C[x_i,x_j]|\geq 3$  and $d(y,C[x_{i+1},x_{j-1}])=0$.

\textbf{Proof}. The proof is by contradiction. Without loss of generality, we may assume that  $x_{n-1}y, x_{a+1}y\in D$, $a\geq 1$ and $d(y,C[x_{1},x_{a}])=0$. Then $a\in [1,n-4]$ (by Lemma 7(i)) and $yx_{a+2}\in D$ (by Claim 3(iii)). From this it is easy to see that 
$$
x_ix_{i-1}\notin D \quad \hbox{for all} \quad  i\not=a+2. \eqno (10) 
$$
We will distinguish two cases, according as $a\geq 2$ or $a=1$.

\textbf{Case 4.1}. $a\geq 2$. Note that $d(y)\leq n-a-1$ (by Claim 1). We extend the path $P_0:=C[x_{a+1},x_{n-1}]$ with vertices $x_{1},x_{2},\ldots , x_{a}$ as much as possible. Then some vertices $z_1,z_2,\ldots , z_d\in \{x_{1},x_{2},\ldots, x_{a}\}$, $d\in [1,a]$, are not on the obtained extended path $P_e$. Using Lemma 2(i), we obtain that $d(z_i)\leq n+d-2$. Therefore 
$$
d(y)+d(z_i)\leq 2n-3+d-a\leq 2n-3. \eqno (11) $$
This together with Lemma 5 implies that $d=1$. Let $z_1:=x_k$. Then $d(x_k)\leq n-1$. First we prove the following Propositions 1 and 2 below\\

\textbf{Proposition 1}. If $x_i\not=x_k$ with $i\in [1,a]$, then $d(x_i)\geq n+a$ and $x_i$ has a partner on $C[x_{a+2},x_{n-1}]$ (i.e., $x_i$ can be inserted into $C[x_{a+2},x_{n-1}])$.

Indeed, the inequality $ d(y)+d(x_k)\leq 2n-2-a$ (by (11) and $d=1$) together with Lemma 5 implies that $d(y)+d(x_i)\geq 2n-1$. Therefore $d(x_i)\geq n+a$, since $d(y)\geq n-a-1$. It is easy to see that 
$$
n+a\leq d(x_i)=d(x_i,C[x_{a+2},x_{n-1}])+d(x_i, C[x_{1},x_{a+1}])\leq d(x_i,C[x_{a+2},x_{n-1}])+2a. 
$$
Hence $d(x_i,C[x_{a+2},x_{n-1}])\geq n-a\geq |C[x_{a+2},x_{n-1}])|+2$, and by Lemma 2(i) the vertex $x_i$ has a partner on $C[x_{a+2},x_{n-1}]$. \fbox \\\\

\textbf{Proposition 2}. Any two vertices $x_i$ and $x_j$ with $k\leq i<j-1\leq a$ (or $1\leq i<j-1\leq k-1$) are non-adjacent.

Indeed, using Proposition 1 and Multi-Insertion Lemma, we obtain that there is an $(x_j,x_i)$-Hamilton ian path, say  $P$,
 and there is an $(x_j,x_i)$-path, say $Q$, 
with vertex set $V(D)-\{x_{j-1}\}$. If $x_jx_i\in D$, then $P$ together with the arc $x_jx_i$ forms a Hamiltonian bypass, and if  $x_ix_j\in D$, then $Q$ together with the arc $x_ix_j$ forms a good cycle, since $d(x_{j-1})\geq n+a$, which contradicts the supposition that $D$ contains no Hamiltonian bypass and good cycle. Therefore $x_j$ and $x_i$ are non-adjacent. \fbox \\\\

Assume first that $k=1$ (i.e., $x_k=x_1$). From Proposition 2 and (10) it follows that
$$
d^-(x_1, C[x_2,x_{a+1}])= d^+(x_a, C[x_1,x_{a-1}])=0. \eqno (12)
$$
In particular, $x_ax_1\notin D$. Thus the triple of vertices $x_1, y, x_a$ satisfies condition $A_0$. Using (11), $d=1$, $z_1=x_1$ and (12), we obtain  
$$
3n-2\leq d(x_1)+d(y)+d^-(x_1)+d^+(x_a)\leq 2n-2-a+d^-(x_1)+d^+(x_a),
$$
and hence 
$$
n+a\leq d^-(x_1)+d^+(x_a)=d^-(x_1, C[x_{a+1},x_{n-1}])+ d^+(x_a, C[x_{a+1},x_{n-1}]).
$$
Now, by Lemma 4, we can insert the path $x_1x_2\ldots x_a$ into $C[x_{a+1},x_{n-1}]$ and obtain an $(x_{a+1},x_{n-1})$-path, say $R$, with vertex set $V(C)$. Therefore, $[x_{a+1}y; Ry]$ is a Hamiltonian bypass, a contradiction.

Assume second that $k\geq 2$ (i.e., $x_k\in C[x_2,x_a]$). From Proposition 2 and (10) it follows that 
$$
d^-(x_1, C[x_2,x_{k}])=0 \quad \hbox{and if} \quad k\leq a-1, \quad \hbox{then} \quad d^+(x_{k}, C[x_{k+2},x_{a+1}])=0, \eqno (13)
$$
$$
d^-(x_1, C[x_2,x_{a+1}])\leq a-k+1 \quad \hbox{and} \quad  d^+(x_{k}, C[x_{1},x_{a+1}])= 1. \eqno (14)
$$
In particular, $x_kx_1\notin D$. The triple of vertices $y,x_k,x_1$ satisfies the condition $A_0$. Hence, using (11), (13) and (14), we obtain  
$$
3n-2\leq d(x_k)+d(y)+d^-(x_1)+d^+(x_k)\leq 2n-2-a+d^-(x_1)+d^+(x_k),
$$
$$
n+a\leq d^-(x_1)+d^+(x_k)= d^-(x_1, C[x_{a+2},x_{n-1}])+ d^+(x_k, C[x_{a+2},x_{n-1}])+$$ $$d^-(x_1, C[x_{1},x_{a+1}]) 
+ d^+(x_k, C[x_{1},x_{a+1}]).
$$
and
$$d^-(x_1, C[x_{a+2},x_{n-1}])+ d^+(x_k, C[x_{a+2},x_{n-1}])\geq n+k-2\geq n. $$
Therefore, by Lemma 4, the path $x_1x_2\ldots x_k$ can be inserted into $C[x_{a+2},x_{n-1}]$. On the other hand, since every vertex $x_i$ with $i\in [k+1,a]$ has a partner on $C[x_{a+2},x_{n-1}]$ (Proposition 1) by Multi-Insertion Lemma there exists an $(x_{a+2},x_{n-1})$-path, say $R$, with vertex set $V(C)-\{x_{a+1}\}$. Therefore, $[x_{a+1}y; x_{a+1}Ry]$ is a Hamiltonian bypass in $D$, which contradicts the supposition that $D$ has no Hamiltonian bypass.\\ 

\textbf{Case 4.2}. $a=1$. Then $x_1$ cannot be inserted into $C[x_{2},x_{n-1}]$. Therefore by Lemma 2(i), $d(x_1)\leq n-1$, and hence
$$
d(y)+d(x_1)\leq 2n-3. \eqno(15)
$$
Recall that $x_2x_1\notin D$ and $x_1x_{n-1}\notin D$ (by (10)). The triples of vertices $y,x_1,x_{n-1}$ and $y,x_1,x_{2}$ satisfies condition $A_0$. Condition $A_0$ together with (15) implies that
$$
3n-2\leq d(x_1)+d(y)+d^+(x_1)+d^-(x_{n-1})\leq 2n-3+d^+(x_1)+d^-(x_{n-1}),
$$
and so
$ d^+(x_1)+d^-(x_{n-1})\geq n+1$.  
A similar argument gives 
$ d^-(x_1)+d^+(x_{2})\geq n+1$.  

The last two inequalities and $d(x_1)\leq n-1$ imply that 
$$ d^-(x_{n-1})+d^+(x_{2})\geq 2n+2-d(x_1)\geq n+3. \eqno (16) 
$$
From $yx_{n-1}\notin D$ (Claim 1), and (10) we obtain that  $d^-(x_{n-1},\{y,x_1,x_2\})\leq 1$ and $d^+(x_{2},\{y,x_1,x_{n-1}\})$ $\leq 2$. This together with  (16) implies that $n\geq 8$ and 
$$
d^-(x_{n-1}, C[x_{3},x_{n-2}])+ d^+(x_2, C[x_{3},x_{n-2}])\geq n > |C[x_{3},x_{n-2}]|+2. $$
By Lemma 4, we can insert the path $x_{n-1}x_1x_2$ into $C[x_{3},x_{n-2}]$ and will obtain an $(x_3,x_{n-2})$-path, say $P$, with vertex set $V(C)$. If $yx_{n-2}\in D$, then $[yx_{n-2}; yP]$ is a Hamiltonian bypass, a contradiction. So, by Lemma 7(i) we can assume that $x_{n-2}$ and $y$ are non-adjacent. From Claim 3(iii) it follows that there exists an integer $b\geq 1$ such that $yx_{n-2-b}\in D$ and $d(y,C[x_{n-b-1},x_{n-2}])=0$. Hence, by Claim 1,
$$ 
 d(y)\leq n-2-b\quad \hbox{and} \quad   d(y)+d(x_1)\leq 2n-(b+3). \eqno (17) $$
It is clear that $n-b-2\not=4$ (Lemma 7(i)). 

Let $n-b-2\geq 5$. Then Claim 3(iii) implies  that $x_{n-b-3}y\in D$. From (17) and Lemma 5 it follows that for every vertex $x_i$ with $i\in [n-b-1,n-2]$ the following inequalities hold
$$
d(y)+d(x_i)\geq 2n-2+(b+3)/2 \quad \hbox{and} \quad d(x_i)\geq n+(3b+3)/2, $$
and hence, using  (10) we obtain 
$$
d(x_i,C[x_3,x_{n-b-3}])\geq n+(3b+3)/2-(2b+4)\geq n-(b+5)/2\geq |C[x_3,x_{n-b-3}]|+2=n-b-3.
$$
Therefore, by Lemma 2(i), every vertex $x_i$, $i\in [n-b-1,n-2]$ has a partner on $C[x_3,x_{n-b-3}]$. By Multi-Insertion Lemma there exists an $(x_3,x_{n-b-3})$-path, say $R$, with vertex set $C[x_3,x_{n-b-3}]\cup C[x_{n-b-1},x_{n-2}]$. Note that $|R|=n-5$. From (16) we have
$$n+3\leq d^-(x_{n-1})+d^+(x_{2})= d^-(x_{n-1},R)+d^+(x_{2},R)+ d^-(x_{n-1}, \{x_1,x_2,y,x_{n-b-2}\})+$$ 
$$d^+(x_{2},\{x_{n-1},x_1,y,x_{n-b-2}\}),
 $$
and, since $d^-(x_{n-1}, \{x_1,x_2,y,x_{n-b-2}\})\leq 2$ and $d^+(x_{2},\{x_{n-1},x_1,y,x_{n-b-2}\})\leq 3$,

$$d^-(x_{n-1},R)+d^+(x_{2},R)\geq n-2\geq |R|+2.$$
By Lemma 4 this means that we can insert the path $x_{n-1}x_1x_2$ into $R$. Therefore there is an $(x_3,x_{n-b-3})$-path, say $Q$, with vertex set $V(C)-\{x_{n-b-2}\}$ and hence , $[yx_{n-b-2}; yRx_{n-b-2}]$ is a Hamiltonian bypass, a contradiction.

Let finally $n-b-2=3$. Then $d(y)=3$, $d(x_1)\leq n-1$, $d^-(x_1)\leq n-3$ and $d^+(x_2)\leq n-2$. Therefore, since $x_2x_1\notin D$, by condition $A_0$ we obtain that  
$$
3n-2\leq d(y)+d(x_1)+d^-(x_{1})+d^+(x_{2})\leq 3n-3,  $$
which is a contradiction, and completes the proof of Claim 4. \fbox \\\\

\textbf{Claim 5}. Let $x_ry, yx_k\in D$ and  $d(y,C[x_{r+1},x_{k-1}])=0$ for some $r, k\in [1,n-1]$, where $3\leq |C[x_{r},x_{k}]|\leq n-2$. Then the vertices $y$ and $x_{k+1}$ are non-adjacent.

\textbf{Proof}. Assume, without loss of generality, that $x_{n-1}y, yx_{a+1}\in D$ (i.e., $x_r=x_{n-1}$ and $x_k=x_{a+1}$) and  $d(y,C[x_{1},x_{a}])=0$ where $a\in [1,n-4]$.

 Suppose that Claim 5 is not true, i.e., the vertices $y$ and $x_{a+2}$ are adjacent. From Lemma 7(i) it follows that $x_{a+2}y\in D$ and $a+2\leq n-3$. Together with Claim 3(iii) this implies that $yx_{a+3}\in D$. It is easy to see that 
$$
x_ix_{i-1}\notin D \quad \hbox{for all} \quad  i\not=a+3. \eqno (18) 
$$
 By Claim 3(iii) there exists a vertex $x_j\in C[x_{1},x_{a}]$ such that $d(x_j)\leq n-1$. Therefore
$$
 d(y)+d(x_j)\leq 2n-(a+2). \eqno (19) $$

\textbf{Proposition 3}. Let $x_l\not=x_j$ with $i\in [1,a]$ (if $a\geq 2$) be an arbitrary vertex. Then $x_l$ has a partner on $C[x_{a+3},x_{n-1}]$ and $d(x_l)\geq n+3a/2$. 

Indeed, by Lemma 5 and (19) the following hold 
$$d(y)+d(x_l)\geq 2n-2 +(a+2)/2 \quad \hbox{and} \quad d(x_l)\geq n+3a/2.$$
Hence, since $x_{l+1}x_l\notin D$ (by (18)), we have that 
$$
n+3a/2\leq d(x_l)=d(x_l, C[x_{a+3},x_{n-1}])+d(x_l, C[x_{1},x_{a+2}])\leq  
d(x_l, C[x_{a+3},x_{n-1}])+2a+1.  $$ 
Therefore
$$
d(x_l, C[x_{a+3},x_{n-1}])\geq n-a/2-1\geq |C[x_{a+3},x_{n-1}])|+2=n-a-1, 
$$  
and by Lemma 2(i), $x_l$ has a partner on $C[x_{a+3},x_{n-1}]$. \fbox \\

Now using Proposition 3, (18) and Multi-Insertion Lemma it is not difficult to show that 
$$
d^+(x_{a+1}, C[x_{j},x_{a}])=d^-(x_{j}, C[x_{1},x_{j-2}]\cup  C[x_{j+1},x_{a+1}])=0, 
$$
(here if $x_j=x_1$ or $x_2$, then $C[x_1,x_{j-2}]=\emptyset$) for otherwise by (18) $a\geq 2$ and $D$ contains 
 a Hamiltonian bypass or a good cycle. In
particular, these equalities imply that 
$$
d^-(x_{j}, C[x_{1},x_{a+2}])\leq 2 \quad \hbox{and} \quad d^+(x_{a+1}, C[x_{1},x_{a+2}])\leq j. \eqno (20)
$$
Note that the triple of vertices $y, x_j, x_{a+1}$ satisfies the condition $A_0$, since $x_{a+1}x_j\notin D$ and the vertices  $y,x_j$ are non-adjacent. The condition $A_0$ together with (19) and (20) implies that 
$$
3n-2\leq d(y)+d(x_j)+d^-(x_j)+d^+(x_{a+1}); $$
$$
n+a\leq d^-(x_j)+d^+(x_{a+1})=d^-(x_j, C[x_1,x_{a+2}])+d^+(x_{a+1},C[x_1,x_{a+2}])+$$
$$d^-(x_j, C[x_{a+3},x_{n-1}])+d^+(x_{a+1},C[x_{a+3},x_{n-1}]).
$$
From this and (20) we obtain that 
$$
d^-(x_j, C[x_{a+3},x_{n-1}])+d^+(x_{a+1},C[x_{a+3},x_{n-1}])\geq n+a-2-j\geq |C[x_{a+3},x_{n-1}]|+2.
$$
Therefore, by Lemma 4, the path $x_jx_{j+1}\ldots x_{a+1}$ has a partner on $C[x_{a+3},x_{n-1}]$. This together with Proposition 3 implies that the path $x_1x_2\ldots x_{a+1}$ has a collection of partners on $C[x_{a+3},x_{n-1}]$, and by Multi-Insertion Lemma there is an $(x_{a+3},x_{n-1})$-path, say $R$, so that $V(R)=V(C)-\{x_{a+2}\}$. This means that $[x_{a+2}y; x_{a+2}Ry]$ is a Hamiltonian bypass, a contradiction. Claim 5 is proved.\\

\textbf{Claim 6}. If $x_ly\in D$ and  $d(y,C[x_{l+1},x_{l+a}])=0$, where $a\in [1,n-4]$, then $yx_{l+a+1}\notin D$.

\textbf{Proof}. The proof is by contradiction. Without loss of generality, assume that $x_{n-1}y\in D$, $d(y,$ $ C[x_1,x_{a}])=0$ and $yx_{a+1}\in D$, where $a\in [1,n-4]$. By Claim 5, the vertices $y$ and $x_{a+2}$ are non-adjacent. If we consider the converse digraph of $D$ we obtain that the vertices $x_{n-2}$ and $y$ also are non-adjacent. It follows from Claim 3(iii) that there is an integer $b\geq 1$ such that $d(y, C[x_{a+2}, x_{a+b+1}])=0$ and $x_{a+b+2}y\in D$. Using the fact that $d(y)\geq 3$ , Lemma 7(i) and again Claim 3(iii) we obtain that $a+b+3\leq n-3$, $yx_{a+b+3}\in D$, and hence
$$
x_ix_{i-1}\notin D \quad \hbox{for all} \quad i\not=a+b+3. \eqno (21) $$
Notice that (by Claim 1)
$$
d(y)\leq n-2-a-b. \eqno (22) $$
On the other hand, by Claim 3(ii) there is a vertex $x_k$ with $k\in [1,a]$ such that $d(x_k)\leq n-1$. This together with (22) implies that 
$$
d(y)+d(x_k)\leq 2n-(a+b+3). \eqno (23)$$
Therefore by Lemma 5, (22) and (23) for every vertex $u\in C[x_1,x_a]\cup C[x_{a+2},x_{a+b+1}]-\{x_k\}$ the following hold
$$
d(u)+d(y)\geq 2n-2+(a+b+3)/2;
$$
$$ d(u)\geq 2n-2+(a+b+3)/2 -n+2+a+b=n+3(a+b+1)/2; \eqno (24) $$
and, since (21),
$$
d(u)=d(u, C[x_1,x_{a+b+2}])+d(u, C[x_{a+b+3},x_{n-1}])\leq d(u, C[x_{a+b+3},x_{n-1}])+2(a+b+1)-1; 
$$
$$
d(u, C[x_{a+b+3},x_{n-1}])\geq n+1-(a+b+1)/2\geq |C[x_{a+b+3},x_{n-1}]|+2=n-a-b-1.
$$
Therefore by Lemma 2(i) the vertex $u$ has a partner on $C[x_{a+b+3},x_{n-1}]$. On the other hand, using this, (24) and Multi-Insertion Lemma it is not difficult to show that 
$$
d^-(x_k,C[x_{1},x_{a+1}])\leq 1 \quad \hbox{and} \quad d^+(x_{a+1},C[x_{k},x_{a+b+2}])=1. 
$$
Hence
$$
d^-(x_k,C[x_{1},x_{a+b+2}])\leq b+2 \quad \hbox{and} \quad  d^+(x_{a+1},C[x_{1},x_{a+b+2}])\leq k. \eqno (25)
$$
Since the triple of vertices $y, x_k, x_{a+1}$ satisfies the condition $A_0$, from (23) and (25) it follows that 
$$
3n-2\leq d(y)+d(x_k)+d^-(x_k)+d^+(x_{a+1})\leq 2n-(a+b+3)+d^-(x_k, C[x_{a+b+3},x_{n-1}])+$$
 $$d^+(x_{a+1},C[x_{a+b+3},x_{n-1}])+b+2+k,$$
and since $k\leq a$,
$$d^-(x_k, C[x_{a+b+3},x_{n-1}])+ d^+(x_{a+1},C[x_{a+b+3},x_{n-1}])\geq 3n-2-2n+(a+b+3)-b-2-k=$$ $$n-1+a-k\geq n-1\geq |C[x_{a+b+3},x_{n-1}]|+2 .$$
Therefore by Lemma 4 the path $x_kx_{k+1}\ldots x_{a+1}$ has a partner on $C[x_{a+b+3},x_{n-1}]$. Thus we have shown that the path $x_1x_2\ldots x_{a+b+1}$ has a collection of partners on $C[x_{a+b+3},x_{n-1}]$. From Multi-Insertion Lemma it follows that there exists an $(x_{a+b+3},x_{n-1})$-path, say $R$, with vertex set $V(C)-\{x_{a+b+2}\}$. Hence, $[x_{a+b+2}y; x_{a+b+2}Ry]$ is a Hamiltonian bypass, which is a contradiction and completes the proof of Claim 6. \fbox \\\\

\textbf{Claim 7}. If $yx_l\in D$ and  $d(y,C[x_{l+1},x_{l+a}])=0$ with $a\in [1,n-4]$, then $x_{l+a+1}y\notin D$.

\textbf{Proof}. Suppose that the claim is not true. Without loss of generality, assume that $yx_{n-1}\in D$, $d(y, C[x_1,x_{a}])=0$ and $x_{a+1}y\in D$, where $a\in [1,n-4]$. Notice that $d(y)\leq n-a-1$ by Claim 1. Lemma 7(i) and  Claims 4 and 6 imply that $yx_{a+2}\in D$ and $x_{n-2}y\in D$. From this it is easy to see that $x_ix_{i-1}\notin D$ for all $i\in [1,n-1]$.

First we prove the following.

\textbf{Proposition 4}. If $d(x_j)\geq n+a-1$ with $x_j\in C[x_1,x_a]$, then $x_j$ has a partner on $C[x_{a+2},x_{n-1}]$ and on $C[x_{a+1},x_{n-2}]$. 

\textbf{Proof of Proposition 4}. Since $x_{j+1}x_j \notin D$, it follows that  $d(x_j, C[x_1,x_{a+1}])\leq 2a-1$. Therefore from
$$
 n+a-1\leq d(x_j)=d(x_j, C[x_1,x_{a+1}])+d(x_j,C[x_{a+2},x_{n-1}]) $$
we obtain that 
$$
d(x_j,C[x_{a+2},x_{n-1}])\geq n-a\geq |C[x_{a+2},x_{n-1}]|+2=n-a,$$
and hence, by Lemma 2(i) $x_j$ has a partner on $C[x_{a+2},x_{n-1}]$. A similar discussion holds  for the path $C[x_{a+1},x_{n-2}]$ and so the proposition is proved. \fbox \\\\

Now we will consider the following cases. 

\textbf{Case 7.1}. $a\geq 2$ and $d(x_k)\leq n+a-1$ for some $x_k\in C[x_1,x_{a}]$. Then, since $d(y)\leq n-a-1$,
$$
d(y)+d(x_k)\leq 2n-2. \eqno (26) $$
Let $x_j\not= x_k$ with $j\in [1,a]$ be an arbitrary vertex. From Lemma 5 and (26) it follows that 
$$
d(y)+d(x_j)\geq 2n-1 \quad \hbox{and} \quad d(x_j)\geq n+a. \eqno (27) 
$$
Without loss of generality, we can assume that $k\geq 2$ (otherwise we consider the converse digraph of $D$). From Proposition 4 it follows that 
$$
d^-(x_{n-1},C[x_1, x_k])=0  \quad \hbox{and} \quad d^+(x_k,C[x_{n-1},x_a])\leq 1  \eqno (28)
$$
(for otherwise, using Multi-Insertion Lemma, we obtain that $D$ contains a Hamiltonian bypass or a good cycle). The triple of vertices $y, x_k, x_{n-1}$ satisfies the condition $A_0$, since $y, x_k$ are non-adjacent and $x_kx_{n-1}\notin D$. Therefore using (26) and (28) we obtain that 
$$
3n-2\leq d(y)+d(x_k)+d^-(x_{n-1})+d^+(x_{k})\leq 2n-2+1+(a-k+1)+$$ $$d^-(x_{n-1}, C[x_{a+1},x_{n-2}])+ d^+(x_{k},C[x_{a+1},x_{n-2}])
$$
and 
$$
d^-(x_{n-1}, C[x_{a+1},x_{n-2}])+ d^+(x_{k},C[x_{a+1},x_{n-2}])\geq n-2-a+k\geq |C[x_{a+1},x_{n-2}]|+2=n-a.
$$
Hence by Lemma 4 we have that the path $x_{n-1}x_1x_2\ldots x_k$ has a partner on $C[x_{a+1},x_{n-2}]$. This together with (27) and Proposition 4 implies that the path $x_{n-1}x_1x_2\ldots x_a$ has a collection of partners on  $C[x_{a+1},x_{n-2}]$, and hence by Multi-Insertion Lemma there is a $(x_{a+1},x_{n-2})$- path, say $R$, with vertex set $V(C)$. Therefore $[x_{a+1}y; Ry]$ is a Hamiltonian bypass, a contradiction.

\textbf{Case 7.2}. $a\geq 2$ and $d(x_j)\geq n+a$ for all $x_j\in C[x_1,x_{a}]$. By Proposition 4 every vertex $x_j$ with $j\in [1,a]$ has a partner on $C[x_{a+2},x_{n-1}]$ and on $C[x_{a+1},x_{n-2}]$. Therefore by Multi-Insertion Lemma, $x_{n-1}$ (respectively, $x_{a+1}$) has no partner on $C[x_{a+1},x_{n-2}]$ (respectively, on $C[x_{a+2},x_{n-1}]$) because of  $x_{a+1}y$ and $x_{n-2}y\in D$ (respectively, $yx_{a+2}$ and $yx_{n-1}\in D$). By Lemma 2(i) this means that

$$
d(x_{n-1},C[x_{a+1},x_{n-2}])\leq n-a-1 \quad \hbox{and} \quad d(x_{a+1},C[x_{a+2},x_{n-1}])\leq n-a-1. \eqno (29)
 $$
On the other hand, using Proposition 4 and  Multi-Insertion Lemma, one can show that $x_{n-1}, x_{a+1}$ are non-adjacent and$$
d^-(x_{1},C[x_{2},x_{a+1}])= d^+(x_{a+1},C[x_{1},x_{a}])=0, \eqno (30)
 $$
$$
d(x_{n-1},\{x_1,x_2,\ldots , x_a,y\})=d(x_{a+1},\{x_1,x_2,\ldots , x_a,y\})=2, 
 $$
since $D$ contains no Hamiltonian bypass and good cycle. 
The last two equalities  together with (29) gives 

$$
d(x_{n-1})\leq n-a+1 \quad \hbox{and} \quad d(x_{a+1})\leq n-a+1. \eqno (31) 
$$
Now using the condition $A_0$, (30) and (31) we obtain that 
$$
3n-2\leq d(x_{n-1})+d(x_{a+1})+d^-(x_{1})+d^+(x_{a+1})\leq 2n-2a+2+d^-(x_{1})+ d^+(x_{a+1})
$$
and 
$$n+2a-4\leq d^-(x_{1})+d^+(x_{a+1})=d^-(x_1, C[x_{a+2},x_{n-1}])+d^+(x_{a+1}, C[x_{a+2},x_{n-1}])+$$
 $$
d^-(x_1, C[x_{2},x_{a+1}])+d^+(x_{a+1}, C[x_{1},x_{a}]\cup \{y\}),
$$
$$
d^-(x_1, C[x_{a+2},x_{n-1}])+d^+(x_{a+1}, C[x_{a+2},x_{n-1}])\geq n+2a-5\geq |C[x_{a+2},x_{n-1}]|+2=n-a,$$
since $a\geq 2$. By Lemma 4 the path $x_1x_2\ldots x_{a+1}$ has a partner on $C[x_{a+2},x_{n-1}]$. Therefore there is an $(x_{a+2},x_{n-1})$-path, say $R$, with vertex set $V(C)$. So we have that $[yx_{n-1}; yR]$ is a Hamiltonian bypass, which is contradiction. This contradiction completes the discussion of the case $a\geq 2$.

\textbf{Case 7.3}. $a=1$. It is easy to see that the arc $x_1x_2$ has no partner on $C[x_3,x_{n-1}]$. Applying Lemma 4 to the arc $x_1x_2$ and to the path $C[x_3,x_{n-1}]$ we obtain that 
$$
d^-(x_1)+d^+(x_2)=d^-(x_1,C[x_3,x_{n-1}])+d^+(x_2,C[x_3,x_{n-1}])+d^-(x_1,\{y,x_2\})+$$ 
$$d^+(x_2,\{y,x_1\})\leq n-1, \eqno (32)
$$
since $d^-(x_1,\{y,x_2\})=0$ and $d^+(x_2,\{y,x_1\})=1$. 
Note that the triple of vertices $y, x_1, x_2$ satisfies condition $A_0$ since $x_1,y$ are non-adjacent and $x_2x_1\notin D$. This together with $d(y)\leq n-2$ and (32) implies that 
$$
3n-2\leq d(y)+d(x_{1})+d^-(x_{1})+d^+(x_{2})\leq d(x_1)+2n-3
$$
and $d(x_1)\geq n+1$. Now by Proposition 4, $x_1$ has a partner on $C[x_3,x_{n-1}]$ and $C[x_2,x_{n-2}]$. Therefore by Multi-Insertion Lemma $x_2$ (respectively, $x_{n-1}$) has no partner on $C[x_3,x_{n-1}]$ (respectively, $C[x_2,x_{n-2}]$). This means that (by Lemma 2(i))

 $$ d(x_2)=d(x_2,C[x_3,x_{n-2}])+d(x_2,\{x_{n-1},x_1,y\})\leq n-1$$ and 
$$
d(x_{n-1})=d(x_{n-1},C[x_3,x_{n-2}])+d(x_{n-1},\{x_{1},x_2,y\})\leq n-1,
$$ 
since $x_{n-1},x_2$ are non-adjacent, $x_1x_{n-1}\notin D$ and $x_2x_1\notin D$. Now using condition $A_0$, (32) and the last two inequalities we obtain  
$$
3n-2\leq d(x_{n-1})+d(x_{2})+d^-(x_{1})+d^+(x_{2})\leq 3n-3,$$
which is a contradiction. Claim 7 is proved. \fbox \\\\

We are now ready to complete the proof of Theorem 12.

From Claims 1 and 2 it follows that there are two distinct vertices $x_k, x_l$ such that $|C[x_k,x_l]|\geq 3$, $y$ is adjacent with $x_k, x_l$ and $d(y,C[x_{k+1},x_{l-1}])=0$. Therefore one of the following cases holds: (i) $x_ky, x_ly\in D$; (ii) $yx_k, yx_l\in D$; (iii) $x_ky, yx_l\in D$; (iv) $yx_k, x_ly\in D$. On the other hand, if $D$ has no Hamiltonian bypass, then  Claims 4-7 imply that each of these cases is impossible. Thus we have a contradiction. The proof of Theorem 12 is completes. \fbox \\\\

\section { Concluding remarks }

Each of  Theorems 1-6  imposes a degree condition on all pairs of non-adjacent vertices (or on all vertices). In the following three theorems imposes a degree condition only for some pairs of non-adjacent vertices. In each of the condition (Theorems 13-16) below $D$ is a strongly connected digraph of order $n$.\\

\textbf{Theorem 13} \cite{[2]} (Bang-Jensen, Gutin, H.Li). Suppose that $min\{d(x),d(y)\}\geq n-1$ and  $d(x)+d(y)\geq 2n-1$ for any pair of non-adjacent vertices $x,y$ with a common in-neighbour, then $D$ is Hamiltonian.\\

\textbf{Theorem 14} \cite{[2]} (Bang-Jensen, Gutin, H.Li). Suppose that $min\{d^+(x)+d^-(y),d^-(x)+d^+(y)\}\geq n$ for any pair of non-adjacent vertices $x,y$ with a common out-neighbour or a common in-neighbour, then $D$ is Hamiltonian.\\

\textbf{Theorem 15} \cite{[3]} (Bang-Jensen, Guo, Yeo). Suppose that $d(x)+d(y)\geq 2n-1$ and  $min\{d^+(x)+d^-(y),d^-(x)+d^+(y)\}\geq n-1$ for any pair of non-adjacent vertices $x,y$ with a common out-neighbour or a common in-neighbour, then $D$ is Hamiltonian.\\

Note that Theorem 15 generalizes Theorem 14. \\
 
In \cite{[10]} the following results were proved: 

(i) if the minimum semi-degree of $D$ at least two and $D$ satisfies the condition of Theorem 13 or 

(ii) $D$ is not directed cycle and satisfies the condition of Theorem 14, then either $D$ contains a pre-Hamiltonian  cycle or $n$ is even and $D$ is isomorphic to the complete bipartite digraph or to the complete bipartite digraph minus one arc with partite sets of cardinalities $n/2$ and $n/2$. 

In \cite{[11]} proved that if $D$ is not directed cycle and satisfies the condition of Theorem 15, then $D$ contains a pre-Hamiltonian  cycle or a cycle of length  $n-2$.\\

We pose the following problem: 

\textbf{Problem}. Characterize those  digraphs which satisfy the condition of Theorem 13 (or 14 or 15) but has no Hamiltonian bypass.

In  \cite{[12]}  the following theorem was proved:

\textbf{Theorem 16}.  Suppose that $min\{d(x),d(y)\}\geq n-1$ and  $d(x)+d(y)\geq 2n-1$ for any pair of non-adjacent vertices $x,y$ with a common in-neighbour. If $n\geq 6$ and the minimum out-degree of $D$ at least two and the minimum in-degree of $D$ at least three, then $D$ contains a Hamiltonian bypass.\\

We believe that Theorem 16 also is true if we require that minimum in-degree at least two instead of three.

\end{document}